\documentclass[twoside]{article}

\usepackage{a4wide, amsmath,  amssymb, amscd, mathptmx, amsthm}
\usepackage{amsfonts}
\usepackage{stmaryrd}
\usepackage{bbm}

\usepackage{color}

\usepackage{pstricks}
\usepackage[all]{xy}\CompileMatrices\SelectTips{cm}{12}

\theoremstyle{plain}
\newtheorem{Thm}{\sc Theorem}[section]
\newtheorem{theorem}[Thm]{\sc Theorem}
\newtheorem{corollary}[Thm]{\sc Corollary}
\newtheorem*{corollary*}{\sc Corollary}
\newtheorem{proposition}[Thm]{\sc Proposition}
\newtheorem*{proposition*}{\sc Proposition}
\newtheorem{lemma}[Thm]{\sc Lemma}

\newtheorem{definition}[Thm]{\sc Definition}

\theoremstyle{remark}
\newtheorem{remark}[Thm]{Remark}
\newtheorem{example}[Thm]{Example}
\newtheorem*{example*}{Example}
\newtheorem*{remark*}{Remark}

\newcommand{\cA}{{\mathcal A}}

\newcommand{\cD}{{\mathcal D}}

\newcommand{\cH}{{\mathcal H}}

\newcommand{\cO}{{\mathcal O}}
\newcommand{\cP}{{\mathcal P}}

\newcommand{\CC}{{\mathbb C}}

\newcommand{\PP}{{\mathbb P}}
\newcommand{\QQ}{{\mathbb Q}}

\newcommand{\VV}{{\mathbb V}}

\newcommand{\ZZ}{{\mathbb Z}}

\newcommand{\bone}{{\mathbbm 1}}

\newcommand{\Mod}[1]{\mathop{{#1}\!\mbox{-}\mathrm{Mod}}}

\newcommand{\End}{{\mathop{ End}\,}}

\newcommand{\cEnd}{{\mathop{\mathcal{E}nd}\,}}

\newcommand{\et}{\mathop{\rm \acute{e}t}}

\newcommand{\cris}{{\mathop{\rm cris}}}

\newcommand{\Spec}{\mathop{\rm Spec \, }}

\newcommand{\gr}{\mathop{{\rm gr}}}

\newcommand{\chr}{\mathop{\rm char}}

\newcommand{\tors}{\mathop{{\rm Torsion}}}
\newcommand{\Hom}{{\mathop{{\rm Hom}}}}
\newcommand{\cHom}{{\mathop{{\cal H}om}}}

\newcommand{\orb}{\mathop{{\rm orb}}}

\newcommand{\GL}{\mathop{\rm GL}}

\newcommand{\Pic}{{\mathop{\rm Pic\, }}}

\newcommand{\topo}{\mathop{\rm top}}

\begin{document}

\markboth{\rm }{\rm  }

\title{On smooth projective D-affine varieties}
\author{Adrian Langer}

\date{\today}

\maketitle

{\sc Address:}\\
Institute of Mathematics, University of Warsaw,
ul.\ Banacha 2, 02-097 Warszawa, Poland\\
e-mail: {\tt alan@mimuw.edu.pl}

\medskip

\begin{abstract} 
{ 
We show various properties of smooth projective D-affine varieties. In particular, 
any smooth projective  D-affine variety is algebraically simply connected and its  image  under a fibration 
is D-affine. In characteristic zero such  D-affine varieties are also uniruled. 

We also show that (apart from a few small characteristics) a smooth projective surface is
D-affine  if and only if it is isomorphic to either $\PP^2$ or $\PP^1\times \PP^1$. In positive 
characteristic, a basic tool in the proof is  a new generalization of Miyaoka's generic semipositivity 
theorem.} 
\end{abstract}

\section*{Introduction}

Let $X$ be a scheme defined over some algebraically closed  field $k$. Let $\cD_X$ be the sheaf of $k$-linear differential operators on $X$. A \emph{$\cD_X$-module} is a left $\cD_X$-module, which is quasi-coherent as an $\cO_X$-module. $X$ is called \emph{D-quasi-affine} if every $\cD_X$-module $M$ is generated over $\cD_X$
by its global sections. $X$ is called \emph{D-affine} if it is D-quasi-affine  and for every $\cD_X$-module $M$ we have $H^i(X, M)=0$ for all $i>0$. 

In \cite{BB} Beilinson and Bernstein proved that every flag variety (i.e., a quotient of a reductive group by some parabolic subgroup) in characteristic $0$ is D-affine.
This fails in positive characteristic (see \cite{KL}), although some flag varieties are still D-affine
(see, e.g., \cite{Ha}, \cite{La} and \cite{Sa}). However, there are no known examples of smooth projective varieties that are D-affine and that are not flag varieties. In \cite{Th} Thomsen proved that any smooth projective toric variety that is D-affine is a product of projective spaces. 

Note that $\cO_X$ has a canonical structure of a $\cD_X$-module coming from the inclusion $\cD_X\subset \End _{k}{\cO_X}$. In particular, if $X$ is a D-affine variety then $H^i(X, \cO_X)=0$ for all $i>0$. This shows that a smooth projective curve  is D-affine if and only if it is
isomorphic to $\PP^1$. However, in higher dimensions this restriction is essentially the only known condition that must be satisfied by D-affine varieties. In the first part of this note we show some 
other properties of smooth projective D-affine varieties. In particular, we prove the following theorem:

\begin{theorem} \label{simply-connected}
Let $X$ be a smooth projective variety defined over an algebraically closed field $k$. 
Let us assume that $X$ is D-affine.  Then the following conditions are satisfied:
\begin{enumerate} 
\item $\pi_1^{\et}(X)=0$.
\item All left $\cD_X$-modules, which are 
coherent as $\cO_X$-modules, are direct sums of finitely many $\cD_X$-modules 
isomorphic to the canonical $\cD_X$-module $\cO_X$.
\item $X$ does not admit any dominant rational map to a curve of genus $\ge 1$.
\item { If $\chr k=0$ then  $X$ is uniruled.}
\end{enumerate}
\end{theorem}

Proof of parts 1 and 2 of Theorem \ref{simply-connected} is divided into two cases depending on the characteristic of the base field. In case of characteristic zero the theorem follows from Theorem \ref{fund-group-char-0}. The proof depends on reducing to the study of unitary representations of the topological fundamental group of $X$. In positive characteristic Theorem \ref{simply-connected} follows from Theorem \ref{fund-group-char-p}. Here we use interpretation of $\cD_X$-modules as stratified bundles. Part 3  follows from part 1 and Proposition \ref{no-maps-to-curves}. 
The last part of the theorem is an application of Miyaoka's theorem  \cite[Corollary 8.6]{Miy} on generic semipositivity of the cotangent bundle of a non-uniruled variety (see Proposition \ref{uniruled} and Remark \ref{vanishing-H1}). In fact, this part of Theorem \ref{simply-connected}  suggests that $X$ should be rationally connected. This problem is studied in Subsection \ref{MRC}, but here we obtain only a partial result on the maximal rationally connected fibration of a D-affine variety.

\medskip

Our next aim is the study of morphisms from smooth D-affine varieties. Here we obtain the following results:

\begin{theorem} \label{main-fibrations}
Let $X$ be a smooth complete variety defined over an algebraically closed field $k$. 
Let us assume that $X$ is D-affine and let $f:X \to Y$ be a surjective morphism to 
some projective variety $Y$.  Then the following conditions are satisfied:
\begin{enumerate} 
\item If $f$ is a fibration then $Y$ does not admit any divisorial contractions.
\item For any $\cD_X$-module $M$ we have $R^if_*M=0$ for $i>0$.
\item If $Y$ is smooth and $f_*\cO_X$ is locally free then for any $\cD_Y$-module $N$ we have $L^if^*N=0$ for $i>0$. 
\item If $f$ is a fibration and $Y$ is smooth  then $Y$ is D-affine.
\end{enumerate}
\end{theorem}

Part 1  follows from the fact that for any effective divisor $D$ the module $\cO_X(*D)$ carries a natural $\cD_X$-module structure  (see Lemma \ref{contractions}).  In case $Y$ is smooth and $k$ has characteristic zero the second part follows from  \cite[2.14 Proposition]{HP} . In general, we use a similar proof following \cite[Theorem 1.4.1]{Ka} (see Proposition \ref{vanishing-for-fibers}). Part 3 follows from part 2 and some calculation in derived categories (see Corollary \ref{image-of-D-affine-is-D-affine}). The last part follows from parts 2 and 3 and again can be found in 
 \cite[2.14 Proposition]{HP} in case $k$ has characteristic zero.
 
\medskip

Part 3 of Theorem \ref{simply-connected}  says that the only smooth projective curve, which is dominated by a smooth projective D-affine variety, is $\PP^1$. In the next part of the paper we  prove that, except possibly for some small characteristics, all smooth projective surfaces that are images of smooth projective D-affine varieties are flag varieties:

\begin{theorem}\label{D-affine-surfaces}
Let $X$ be a smooth projective variety defined over an algebraically closed field $k$
and let $f: X\to Y$ be a fibration over a smooth projective surface. If $\chr k=0$ or $\chr k>7$
and $X$ is D-affine then { $f$ is flat and}  $Y=\PP^2$ or $Y=\PP^1\times \PP^1$.
\end{theorem}

Let us recall that products of projective spaces are D-affine in any characteristic  (see \cite[Korollar 3.2]{Ha} or \cite{Th}). In particular, a smooth projective surface defined over  an algebraically closed field  of characteristic $0$ or $>7$ is 
D-affine if and only if it is isomorphic to either  $\PP^2$ or $\PP^1\times \PP^1$.

We prove a slightly more precise result in Section \ref{surfaces}. Let us mention that recently D. Rumynin in \cite{Ru} proved that the only D-affine rational surfaces are flag varieties.

In characteristic $0$ the above classification result follows from Theorem \ref{simply-connected}
and Theorem \ref{main-fibrations}. However, the positive characteristic case is more delicate and we need the following  positive characteristic  version of Miyaoka's generic semipositivity theorem. 

\begin{theorem} \label{simple-sym-pow}
Let $X$ be a smooth projective surface defined over an algebraically closed field 
$k$ of characteristic $p$. Let us fix an ample divisor $H$ on $X$ and assume that $p>K_XH+1$. 
{ If $X$ is not uniruled then $\Omega_X$ is generically $H$-semipositive.}
\end{theorem}

{ 
For the definition and basic properites of generically semipositive sheaves in positive characteristic see
Subsection \ref{gen-semipositive}. The most important property is their good behaviour under various tensor operations like symmetric or divided powers. 

The only known result on generic semipositivity of cotangent bundle for non-uniruled varieties in positive characteristic concerns varieties with trivial canonical divisor (see \cite[Theorem 0.1]{La-GP}).
However, it does not say anything about the most interesting case of varieties of general type. This is covered by the above theorem but only in the surface case. The higher dimensional version seems to require different techniques. A more precise version of Theorem \ref{simple-sym-pow} is contained in Theorem \ref{symmetric-power}. We show that this generalization is optimal (see Subsection \ref{Ekedahl-example} and Remark \ref{Eke-rem}).
}

\medskip

The structure of the paper is as follows. In Section 1 we recall some auxiliary results. 
{  In Section 2 we prove several facts about tensor operations and generic semipositivity of sheaves in positive characteristic.
In Section 3 we prove Theorem \ref{simply-connected}.  In Section 4 we prove Theorem \ref{main-fibrations}.
In Section 5 we study uniruledness of surfaces  in positive characteristic proving a generalization of Theorem \ref{simple-sym-pow}. Finally, in Section 6 we use these results to study smooth projective surfaces that are images of D-affine varieties, proving a generalization of Theorem \ref{D-affine-surfaces}. We also make some remarks on the case of D-affine $3$-folds.}

\medskip

\subsection*{Notation}

Let $X$ and $Y$ be algebraic varieties defined over an algebraically closed field.
 
A \emph{divisorial contraction} is a proper birational morphism $f:X\to Y$, which contracts some divisor $D$ to a subscheme of codimension $\ge 2$ and that is an isomorphism outside of $D$. 

A \emph{fibration} is a morphism $f:X\to Y$ such that $f_*\cO_X=\cO_Y$ (in particular, we allow $f$ to be the identity or a birational morphism).

Let us assume that $X$ is a smooth projective variety  and let us fix an ample divisor $H$ on $X$.
If $E$ is a torsion free coherent $\cO_X$-module then  by $\mu_{\max, H} (E)$  we denote the slope of the maximal destabilizing subsheaf of $E$ (with respect to $H$). Similarly, we use $\mu_{\min, H} (E)$ to denote the slope of the minimal destabilizing quotient of $E$.  When it is clear from the context which polarization is used, we omit $H$ in the notation and write  $\mu_{\max} (E)$ and  $\mu_{\min} (E)$ instead of  $\mu_{\max, H} (E)$ and  $\mu_{\min, H} (E)$, respectively.

\section{Preliminaries}

{ \subsection{$\cA$-affine varieties}

Let $X$ be a smooth variety defined over an algebraically closed field $k$. 
Let $\cA$ be any sheaf of rings on $X$ with a ring homomorphism $\cO_X\to \cA$ such that the image of 
$k\to \cO_X\to \cA$ is contained in the center of $\cA$. Let us also assume that $\cA$ is quasi-coherent as a left $\cO_X$-module. In the following by an \emph{$\cA$-module} we mean a left $\cA$-module, which is quasi-coherent as an $\cO_X$-module.

\begin{definition}
We say that $X$ is \emph{$\cA$-quasi-affine} if  any $\cA$-module is generated over $\cA$ by its global sections. We say that $X$ is \emph{almost $\cA$-affine} if for any $\cA$-module $M$ we have $H^i(X, M)=0$ for all $i>0$.
$X$ is \emph{$\cA$-affine} if it is both $\cA$-quasi-affine and almost $\cA$-affine.
\end{definition}

If $\cA=\cD_X$ we talk about D-quasi-affine, almost D-affine and D-affine varieties, respectively.

Let us recall that if $M$ is an $\cA$-module then $\Gamma (X, M)$ is a $\Gamma (X,A)$-module. This module has an induced  $\Gamma(X, \cO_X)$-module structure, which agrees with the $\Gamma(X, \cO_X)$-module structure on sections of  $M$ considered as an $\cO_X$-module. Similar fact holds for the derived functor. So in the following we can check vanishing of the derived global sections $H^i(X, M)$ of $M$ treated as an $\cO_X$-module (or even as the derived functor of sections of $M$ treated as a sheaf of abelian groups).

\medskip

If $\cA=\cO_X$ then Serre's theorem says that an almost $\cA$-affine variety is also $\cA$-affine. This fails for more general sheaves of rings. For example, by the Beilinson--Bernstein theorem  this fails for certain rings of twisted differential operators on flag varieties (see, e.g., \cite[Lemma 7.7.1]{Ka} for an explicit example).

\medskip

A special case of $\cA$ is that of  the universal enveloping algebra of some Lie algebroid. By definition such 
$\cA$ comes equipped with an $\cO_X$-linear morphism of sheaves of rings $\cA\to \cD_X$.
The following proposition shows that classification of D-affine varieties gives also classification of such
$\cA$-affine varieties:

\begin{proposition}
Assume that there exists a morphism of sheaves of rings $\cA\to \cD_X$, which is compatible with left $\cO_X$-module structures. 
If $X$ is $\cA$-affine then it is also  D-affine.
\end{proposition}

\begin{proof}
Assume that $X$ is $\cA$-affine. If $M$ is a $\cD_X$-module then it has also an induced $\cA$-module structure and hence $H^i(X, M)=0$ for all $i>0$. 
By \cite[Proposition 1.5.2]{Ka} (or \cite[Proposition 1.4.4]{HTT}) a variety $X$ is D-affine if and only if it is almost D-affine and for any non-zero $\cD_X$-module $M$ we have $\Gamma (X, M)\ne 0$.
Thus  it is sufficient to check that for any non-zero $\cD_X$-module $M$ we have $\Gamma (X, M)\ne 0$. But again such $M$ has an induced $\cA$-module structure and by $\cA$-affinity, $M$ is generated over $\cA$ by its global sections. In particular,   $\Gamma (X, M)\ne 0$ as required. 
\end{proof}

A special case when the above proposition applies is when $B$ is a simple normal crossing divisor and 
$\cA$ is the universal enveloping algebra of the Lie algebroid $T_X(-\log B)\subset T_X$. The proposition 
shows that  "log D-affine varieties" are D-affine.

\medskip

Apart from the usual sheaf $\cD_X$ of $k$-linear differential operators one can also consider the sheaf $\cD_X^{\cris}$ of $k$-linear crystalline differential operators on $X$. This is defined as the universal enveloping algebra of the tangent Lie algebroid $T_X$.  
There exists a canonical morphism $\cD_X^{\cris}\to \cD_X$ of sheaves of rings.
If $\chr k=0$ then this morphism is an isomorphism. However, if $\chr k>0$ then this morphism is neither injective nor surjective. In this case the basic difference between $\cD_X$ and $\cD_X^{\cris}$ is that whereas for $\cD_X$ the sheaf associated to the standard order filtration is isomorphic to $\bigoplus (S^iT_X^*)^*$,  for $\cD_X^{\cris}$ 
the sheaf associated to the standard order filtration is isomorphic to $\bigoplus S^iT_X$.

The following proposition shows that  $\cD_X^{\cris}$-affinity in positive characteristic is a trivial notion.

\begin{proposition}
Let $X$ be a smooth projective variety defined over an algebraically closed field $k$ of positive characteristic $p$. 
If $X$ is $\cD_X^{\cris}$-affine then $X$ is a point.
\end{proposition}

\begin{proof}
Assume that $\dim X>0$ and let $L$ be a very ample line bundle on $X$. Then the Frobenius pull-back $M=F_X^*L^{-1}$ carries a canonical integrable connection, giving $M$ a left $\cD_X^{\cris}$-module structure. 
Since $X$ is $\cD_X^{\cris}$-affine we have $\Gamma (X, M)\ne 0$. But $M=L^{-p}$ and $\Gamma (X, L^{-p})=0$, a contradiction.
\end{proof}

}

{ \subsection{D-affinity}}

We will often use the fact that if $X$ is D-affine and $M\ne 0$ is a $\cD_X$-module then 
$\Gamma (X, M)\ne 0$. This follows immediately from the definition of a D-quasi-affine variety.
In fact, we have the following more general proposition (see, e.g., \cite[Proposition 1.4.4]{HTT}):

\begin{proposition} \label{equivalence}
Let $X$ be a D-affine variety defined over some algebraically closed field. Then the functor
$$\Gamma (X, \bullet ): \Mod{ \cD_X}\to \Mod{ \Gamma (X, \cD_X)}$$
is an equivalence of categories with a quasi-inverse given by 
$$ \cD_X \otimes _{\Gamma (X, \cD_X)} \bullet: \Mod{\Gamma (X, \cD_X)} \to \Mod{ \cD_X} .$$
\end{proposition}

The following lemma is well-known (see \cite[Proposition  2.3.3]{Ha}), but we recall { its} proof 
for the convenience of the reader. It is an analogue of the fact that a quasi-affine variety $X$ is affine if and only if 
$H^i(X, \cO_X)=0$ for all $i>0$.

\begin{lemma}\label{Serre}
Let $X$ be D-quasi-affine. Then $X$ is D-affine if and only if  $H^i(X, \cD_X)=0$
for all $i>0$.  
\end{lemma}

\begin{proof}
By Grothendieck's vanishing theorem for every $\cD_X$-module $M$ we have $H^i(X, M)=0$ for $i$
larger than the dimension of $X$. So it is sufficient to prove that for $n\ge 1$ if for all $\cD_X$-modules $M$  we have $H^i(X, M)=0$ for $i>n$ then for all $\cD_X$-modules $M$ 
we have $H^i(X, M)=0$ for $i\ge n$. Since a  $\cD_X$-module $M$ is globally generated as a $\cD_X$-module
we have a short exact sequence
$$0\to N\to \cD_X\otimes_{\Gamma (X, \cD_X)}{\Gamma (X, M)}\to M\to 0$$
for some  $\cD_X$-module $N$. From the long exact cohomology sequence we have
$$ H^i(X,\cD_X) \otimes_{\Gamma (X, \cD_X)}{\Gamma (X, M)}\to H^i(X,M)\to H^{i+1}(X,N),$$
which proves the required implication.
\end{proof}

{

The following lemma is a small generalization of \cite[Lemma 1]{Th}. 

\begin{lemma}\label{Thomsen}
Let $X$ be a smooth variety defined over an algebraically closed field $k$. Let 
$U$ be an open subset of $X$ such that its complement in $X$ is non-empty and has codimension $1$.
Let $j: U\hookrightarrow X$ be the corresponding embedding.
Assume that $X$ is D-affine.
Then the restriction map  $j^*: \Gamma(X, \cO_X) \to \Gamma(U, \cO_U)$ is not an isomorphism. Moreover, if  
$X\backslash U$ has pure codimension $1$ in $X$ then
$$H^i(U, \cO_U)=0$$
for all $i>0$. In particular, $U$ is quasi-affine if and only if it is affine.
\end{lemma}

\begin{proof}
Since $\cO_U$ is a $\cD_U$-module, $j_*\cO_U$ is a $\cD_X$-module (see \cite[Example 1.5.22 and Proposition 1.5.29]{HTT}). Since the canonical map $\cO_X\to j_*\cO_U$ of $\cD_X$-modules is not an isomorphism (even of $\cO_X$-modules), the corresponding map on global sections 
$ \Gamma(X, \cO_X)\to \Gamma (X, j_*\cO_U)=\Gamma(U, \cO_U)$ 
is not an isomorphism (as $\Gamma (X, \bullet):  \Mod{ \cD_X} \to \Mod{\Gamma (X, \cD_X)}$ is an equivalence of categories by Proposition \ref{equivalence}). 

Now if we assume that $X\backslash U$ has pure codimension $1$ in $X$ then $j$ is an affine morphism, so 
$$H^i(U, \cO_U)= H^i(X, j_*\cO_U)=0$$
for all $i>0$. The last part follows from the criterion similar to the one from Lemma \ref{Serre}.
\end{proof}
}

{ \subsection{Simply connected varieties} }

In proof of Theorem \ref{simply-connected} we need the following proposition:

\begin{proposition} \label{no-maps-to-curves}
Let $X$ be a smooth projective variety defined over an algebraically closed field $k$. If $\pi_1^{\et} (X)=0$
then  $X$ does not admit any dominant rational map to a curve of genus $\ge 1$. Moreover, if there exists a fibration $f: X\to \PP ^1$ then it has at most two multiple fibers.
\end{proposition}

\begin{proof}
Let $f: X\dasharrow  C$ be a dominant rational map to a smooth projective curve $C$.
Note  that $f$ extends to a morphism on an open subset $U$ such that the complement of $U$ in 
$X$ has codimension $\ge 2$. This follows from the fact that a rational map from a smooth curve to a projective variety always extends to a morphism.
Then $ \pi_1^{\et}(U)=0$ by \cite[Expos\'e X, Corollaire 3.3]{SGA1}.
Taking normalization of the graph of $f$ we can find a normal projective variety $\tilde X$, a birational morphism $\tilde X\to X$, which is an isomorphism over $U$, and a morphism $\tilde f: \tilde X\to C$.
Note that $\pi_1^{\et}(U)\to \pi_1^{\et}(\tilde X)$ is surjective, so $\tilde X$ is algebraically simply connected. Let us consider the Stein factorization of $\tilde f$
$$\tilde f: \tilde X\mathop{\longrightarrow}^{g} D\mathop{\longrightarrow}^{h}  C.$$
$D$ is  a smooth projective curve and $g_*\cO_X=\cO_D$, i.e., $g$ is a fibration.
Then we have a surjective map $\pi_1^{\et}(\tilde X)\to \pi_1^{\et}(D)$. Therefore $ \pi_1^{\et}(D)=0$ and we get $D=\PP^1$. In particular, $C=\PP^1$ and $h$ is a finite covering.

Now let us assume that $f: X\to C=\PP ^1$ is a fibration and let us consider all the points $Q_i\in \PP^1$ such that $f$ has multiple fibres of multiplicity $m_i$ over $Q_i$. 

If $k=\CC$ then we have a surjective map from $\pi_1^{\topo}(X)$ to the orbifold fundamental group $\pi_1^{\orb}(C_f)$ of $C$ with respect to $f$ (see \cite[Theorem 2.1]{LS}). This last group is defined as the quotient of  $\pi_1^{\topo}(C-\{Q_i\}_i)$ by the normal subgroup generated  by all the elements of the form $\gamma_i^{m_i}$, where $\gamma _i$
is a simple loop going around the point $Q_i$. But then we get a surjective map from $\pi_1^{\et}(X)$
to the profinite completion of $\pi_1^{\orb}(C_f)$. This last group is clearly non-zero if $C=\PP^1$ and $f$ has at least $3$ multiple fibers.

The proof in an arbitrary characteristic is analogous. Here one can define the \'etale orbifold fundamental group 
$\pi_1^{\orb, \et}(C_f)$ and prove that there exists a surjective homomorphism  $\pi_1^{\et}(X)\to \pi_1^{\orb, \et}(C_f)$
(see \cite[Definition 4.25]{Mi}; note however that  by \cite[Remark 2.2]{LS} the sequence from  \cite[Theorem 4.22]{Mi} is non-exact in the non-proper case).  Since  $\pi_1^{\et}(X)=0$, we have $\pi_1^{\orb, \et}(C_f)=0$. 
Again, one shows that this implies that $f$ has at most two multiple fibers (see \cite[Theorem 1.3]{Mi}).
\end{proof}

{
\section{Semistability and generic semipositivity of sheaves}

Let us fix a normal projective variety $X$ defined over an algebraically closed field $k$
and an ample divisor $H$ on $X$. In this section we gather several facts about strong semistability 
of sheaves in positive characteristic.

\subsection{Bounds on semistability of tensor products}

If $\chr k=p$ then we denote by $F_X:X\to X$ the  absolute Frobenius morphism.

Let $E$ be a torsion free coherent $\cO_X$-module. 
Then we define 
$$L_{\max , H} (E):=\left\{ 
\begin{array}{cl}
\lim _{m\to \infty} \frac{\mu_{\max , H} ((F_X^m)^*E)}{p^m}&\quad \hbox{ \rm if $\chr k=p$,}\\
\mu_{\max , H} (E)& \quad \hbox{ \rm if $\chr k=0$.}\\
\end{array}
\right.$$
Similarly, we can define $L_{\min , H} (E)$. 
Both $L_{\max , H} (E)$ and $L_{\min , H} (E)$ are well defined rational numbers (see \cite[2.3]{La0}).
We say that $E$ is \emph{strongly slope $H$-semistable} if $L_{\max , H} (E)=\mu _{\max , H} (E)$.

\medskip

Let $\rho: \GL (r)\to \GL (s)$ be a representation mapping the centre of  $\GL (r)$ to the centre of $\GL (s)$.

If $E$ is a rank $r$  torsion free coherent $\cO_X$-module then its reflexivization $E^{**}$ is locally free on an open subset $j:U\hookrightarrow X$ such that its complement in $X$ has codimension $\ge 2$. Let $\cP$ be a principal 
$\GL(r)$-bundle on $U$ associated to $j^*(E^{**})$ and let $\cP_{\rho}$ be the principal $\GL(s)$-bundle on $U$
obtained from $\cP$ by extension of structure group via $\rho$. We can associate to $\cP_{\rho}$ a rank $s$ locally free $\cO_U$-module $E_{\rho}$. Then we set $\hat E_{\rho}= j_*E_{\rho}$. By definition $\hat E_{\rho}$ is a reflexive sheaf.

In the following we will need the following theorem of Ramanan and Ramanathan (see \cite[Theorem 4.9]{La0}).

\begin{theorem}\label{RR}
If $E$ is strongly slope $H$-semistable then $\hat E_{\rho}$ is also strongly slope $H$-semistable.
\end{theorem}

For two torsion free coherent $\cO_X$-modules  $E_1$ and $E_2$ we denote by $E_1{\hat \otimes} E_2$  the reflexivization of $E_1\otimes E_2$. Similarly, if $E$ is a  torsion free coherent $\cO_X$-module then we set 
$\hat \bigwedge^j E= (\bigwedge^jE)^{**}$, $\hat S^j E= (S^jE)^{**}$  and $E^{\hat \otimes j}=E {\hat \otimes} ...{\hat \otimes} E$, where $E$ appears in the product $j$ times. Note that  the $j$-th divided power ${\Gamma}^j E=(S^jE^*)^*$ is already reflexive so we do not introduce a new notation for its reflexivization.

As a corollary of the above theorem one gets the following result:

\begin{corollary}\label{ss-tensor-product}
\begin{enumerate}
\item If $E$ is strongly slope $H$-semistable then $E^*$, $\cEnd _{\cO_X} E$,  $E^{\hat \otimes j}$, $\hat \bigwedge^j E$,
$\hat S^j E$ and ${\Gamma}^j E$ are also strongly slope $H$-semistable.
\item If $E_1$ and $E_2$ are strongly slope $H$-semistable then $E_1{\hat \otimes} E_2$ is strongly slope $H$-semistable.
\end{enumerate}
\end{corollary}

\begin{proof}
The first part is obtained by applying Theorem \ref{RR} to the corresponding representation, e.g., 
$\cEnd _{\cO_X} E$ is equal to  $\hat E_{\rho}$  for the adjoint representation of $\GL (r)$, and $\hat S^j E$ is equal to  $\hat E_{\rho}$  for the symmetric representation $\GL(r)=\GL (V)\to \GL (S^jV)$.

To prove the second part let us note that if $\det E_1=\det E_2=\cO_X$ then 
 $E_1\hat \otimes E_2$ is strongly slope $H$-semistable as it is a direct summand of $\cEnd _{\cO_X} (E_1\oplus E_2^*)$.
 Now let $r_i$ be the rank of $E_i$ for $i=1,2$.
 If there exist line bundles $L_1$ and $L_2$ such that $\det E_1=L_1^{r_1}$ and $\det E_2=L_2^{r_2}$
 then $E_1\otimes L_1^{-1}$ and  $E_2\otimes L_2^{-1}$ are as in the previous case so their tensor product is 
 strongly slope $H$-semistable. This implies that $E_1\hat \otimes E_2$ is also strongly slope $H$-semistable. 
 
Now let us consider the general case. By the Bloch--Gieseker covering trick (see \cite[Lemma 2.1]{BG})
there exists a normal projective variety $\tilde X$ and a finite flat surjective covering $f: \tilde
X\to X$ together with line bundle $L_1$ and $L_2$ such that $f^*(\det
E_i)^{-1}=L_i^{r_i}$ for $i=1,2$. Then $f^*E_i$ are strongly slope $f^*H$-semistable, so by the above 
$f^*E_1\hat \otimes f^*E_2$ is also strongly slope $f^*H$-semistable. This implies that $E_1\hat \otimes E_2$ is also
strongly slope $H$-semistable.
\end{proof}

The following theorem is a corollary of the Ramanan--Ramanathan theorem (Theorem \ref{RR}) and the author's results 
(see, e.g., \cite[Theorem 2.13]{La0}). Proof of the first part of the theorem was indicated by the author in 
 \cite[2.3.3]{La0}.

\begin{theorem} \label{L_max-tensor-product}
\begin {enumerate}
\item Let $E_1$ and $E_2$ be torsion free coherent $\cO_X$-modules. Then we have
$$L_{\max, H} (E_1{\hat \otimes} E_2)=L_{\max, H}( E_1)+L_{\max ,H} (E_2).$$
\item Let $E$ be a torsion free coherent $\cO_X$-module.
Then 
$$ L_{\max, H} ({ \hat S}^j E)= L_{\max, H} ({ \Gamma}^j E)=  j \, L_{\max, H} (E).$$
\end{enumerate}
Similar equalities hold if we replace $L_{\max}$ by $L_{\min}$.
\end{theorem}

\begin{proof}
By \cite[Theorem 2.13]{La0} for all large $m\ge 0$ the quotients of the Harder--Narasimhan filtrations of 
$(F_X^m)^*E_1$ and $(F_X^m)^*E_2$ are strongly slope $H$-semistable.
Let $F_{\bullet}^m$ be the Harder--Narasimhan filtration of $(F_X^m)^*E_1$ 
and let $G_{\bullet}^m$ be the Harder--Narasimhan filtration of $(F_X^m)^*E_1$ . Then 
$(F_X^m)^*E_1{\hat \otimes} (F_X^m)^*E_2$
has a filtration whose quotients agree with  tensor products $\gr _{F^m}^i ((F_X^m)^*E_1)\hat \otimes \gr _{G^m}^j ((F_X^m)^* E_2)$ outside of a closed subset of codimension $\ge 2$. By Corollary \ref{ss-tensor-product}
all these quotients are strongly $H$-semistable. So we have
$$
\mu _{\max, H} ((F_X^m)^*E_1{\hat \otimes} (F_X^m)^*E_2)\le \max_{i,j}( \mu ( {\gr} _{F^m}^i ((F_X^m)^*E_1)) + \mu ( {\gr} _{G^m}^i ((F_X^m)^*E_1))  ))=p^m(L_{\max, H}( E_1)+L_{\max ,H} (E_2)).
$$
 Since  ${\gr} _{F^m}^1 ((F_X^m)^*E_1) {\hat \otimes} {\gr} _{G^m}^1 ((F_X^m)^*E_2)$
is a subsheaf of $(F_X^m)^*E_1{\hat \otimes} (F_X^m)^*E_2$ we also have inequality 
$$
\mu _{\max, H} ((F_X^m)^*E_1{\hat \otimes} (F_X^m)^*E_2)\ge 
\mu ( {\gr} _{F^m}^1 ((F_X^m)^*E_1)) + \mu ( {\gr} _{G^m}^1 ((F_X^m)^*E_1)))=p^m(L_{\max, H}( E_1)+L_{\max ,H} (E_2)).
$$
Thus we get the first equality. The analogous equality for $L_{\min}$ is proven in an analogous way.

The proof of the second part of the theorem is similar. Let us first consider the case of the symmetric powers.
By \cite[Theorem 2.13]{La0}  for large $m\ge 0$ all the quotients of the Harder--Narasimhan filtration $F_{\bullet}^m$  of  $(F_X^m)^*E$ are strongly slope $H$-semistable. Assume that there are exactly $s$ factors in this filtration.
Then $S^i ((F_X^m)^*E)$ has a filtration with quotients isomorphic outside of a closed subset of codimension $\ge 2$
to
$$\hat S^{i_1}({\gr} _{F^m}^1 ((F_X^m)^*E))\hat\otimes  \hat S^{i_2}({\gr} _{F^m}^2 ((F_X^m)^*E)) \hat \otimes ...\hat \otimes \hat S^{i_s}({\gr} _{F^m}^s ((F_X^m)^*E)),$$
where $i_1+...+i_s=j$. By Corollary \ref{ss-tensor-product} all these quotients are strongly $H$-semistable 
and as before one can easily see that 
$$\mu _{\max, H} (\hat S  ^j(  (F_X^m)^*E  ) )= \mu(\hat S^{j}({\gr} _{F^m}^1 ((F_X^m)^*E)) )= p^m j L_{\max, H} (E).$$
This implies $ L_{\max, H} ({ \hat S}^j E)=  j \, L_{\max, H} (E).$ Equality for $L_{\min}$ is analogous.

Now the equality for divided powers follows from 
$$ L_{\max, H} ({ \Gamma}^j E) =- L_{\min, H} ({ \hat S}^j (E^*))= - j \, L_{\min, H} (E^*)= j L_{\max, H} (E)$$
and the similar equalities for $L_{\min}$.
\end{proof}

\subsection{Generically semipositive sheaves} \label{gen-semipositive}

Let $X$ be a smooth projective variety defined over an algebraically closed field $k$
and let $H$ be a fixed ample divisor on $X$. The following definition comes from \cite[Definition 1.6]{La-GP}.

\begin{definition}\label{def-gen-semipositive}
A  torsion free coherent $\cO_X$-module $E$ is \emph{generically $H$-semipositive} if $L_{\min , H} (E)\ge 0$.
\end{definition}

If $\chr k=0$ then this definition coincides with the usual definition of generically $H$-semipositive sheaves. 
Let us also recall that in positive characteristic it is not known if the restriction of a generically $H$-semipositive sheaf to a general complete intersection curve $C\in |m_1H|\cap ...\cap |m_{n-1}H|$ with $m_i\gg 0$ is still generically semipositive. However, generically semipositive sheaves are still well behaved with respect to tensor operations, etc.
More precisely, generically semipositive sheaves satisfy the following properties:

\begin{proposition}
\begin{enumerate}
\item Let $$0\to E_1\to E\to E_2\to 0$$ be a short exact sequence of torsion free coherent $\cO_X$-modules.
If $E$ is generically $H$-semipositive then $E_2$ is generically $H$-semipositive. If $E_1$ and
 $E_2$ are generically $H$-semipositive then $E$ is generically $H$-semipositive.
\item If $E_1$ and $E_2$ are generically $H$-semipositive then $E_1{\hat \otimes} E_2$ is generically 
$H$-semipositive.
\item If $E$ is generically $H$-semipositive then for all positive integers $j$ the sheaves $E^{\hat \otimes j}$, 
$\hat S^j E$ and ${\Gamma}^j E$ are generically $H$-semipositive.
\end{enumerate}
\end{proposition}

\begin{proof}
The first assertion follows from the fact that 
$$L_{\min, H} (E)\ge \min (L_{\min, H} (E_1), L_{\min, H} (E_2)).$$
The second and third assertion follow directly from Theorem \ref{L_max-tensor-product}.
\end{proof}
}

\section{Proof of Theorem \ref{simply-connected}}

In this section we prove Theorem \ref{simply-connected}.  The  proof is divided into two cases depending on the characteristic of the base field.

\subsection{Theorem \ref{simply-connected} in the characteristic zero case}

In this subsection we prove parts 1, 2 and { 4} of Theorem \ref{simply-connected} in case the base field $k$ has characteristic zero. Part 3 and the second assertion in 4 follow from 1 and Proposition \ref{no-maps-to-curves}. First, let us prove the last part of  Theorem \ref{simply-connected}:

{ 
\begin{proposition}\label{uniruled}
Let $X$ be a smooth projective variety defined over an algebraically closed field $k$ of characteristic $0$.
If $X$ is D-quasi-affine then it is uniruled.
\end{proposition}

\begin{proof}
Let us fix an ample line bundle $L$ on $X$. Since $M=\cD_X\otimes _{\cO_X} L^{-1}$ is a left $\cD_X$-module,
we have $\Gamma (X, M)\ne 0$. Note that $M$ has a natural good filtration by coherent $\cO_Y$-submodules $F_iM:=\cD_X^{\le i}\otimes _{\cO_Y}L^{-1}$, where $\cD_X^{\le i}$ denotes the sheaf of differential operators of order $\le i$.  In particular, there exists some $i\ge 0$ such that 
$\Gamma (X, F_i M)\ne 0$. Since $F_0M \subset F_1M\subset ...\subset F_iM$,
there exists some $j\le i$ such that $\Gamma (X, F_jM/F_{j-1}M)\ne 0$. But $F_jM/F_{j-1}M= S^j T_X\otimes L^{-1})$. Thus there exists some $j>0$ such that $S^jT_X$ contains $L$ as an $\cO_X  $-submodule.

If $X$ is not uniruled then by Miyaoka's theorem \cite[Corollary 8.6]{Miy} $\Omega_X$ is generically semipositive. In other words, for any fixed ample polarization we have $\mu_{\min} (\Omega_X)\ge 0$.
Since tensor operations on semistable sheaves preserve semistability, this inequality implies that
$\mu_{\min} (S^j\Omega_X)\ge 0$ { (see { Theorem \ref{L_max-tensor-product}})}. But then $\mu_{\max} (S^jT_X)\le 0$, which contradicts the fact that
$S^jT_X$ contains an ample line bundle. 
\end{proof}
}

\medskip

Now let us go back to the proof of parts 1 and 2 of Theorem \ref{simply-connected}.  Without loss of generality we can assume that $k=\CC$. 
If $X$ is smooth complex projective variety and $X$ is D-affine then $\cO_X$ has only one structure of a $\cD_X$-module as $h^0(X, \Omega_X)=h^1(X, \cO_X)=0$. However, for a left $\cD_X$-module $M$ the evaluation map
$$\cO_X\otimes _{\CC}\Gamma (X, M)\to M$$
is usually not a map of $\cD_X$-modules (if it were, one could easily see that it is an isomorphism of $\cD_X$-modules, proving that $X$ is D-affine). The idea behind the proof of the following theorem is that if $M$ is a locally free $\cO_X$-module of finite rank underlying a unitary representation then this map is a non-trivial map between slope semistable bundles of degree $0$ (with respect to some polarization) and we get enough information to prove the first part of Theorem \ref{simply-connected}. Over complex numbers, it is easy to see that this implies the second part of Theorem \ref{simply-connected}.

Before giving the proof, let us recall that every $\cD_X$-module, which is coherent as an $\cO_X$-module, 
is locally free as an $\cO_X$-module (see \cite[2.15 and 2.17]{BO} or \cite[Theorem 1.4.10]{HTT}). Moreover, giving 
a left $\cD_X$-module structure extending a given $\cO_X$-modules structure is equivalent 
to giving an integrable connection. So left $\cD_X$-modules, which are coherent as $\cO_X$-modules,
correspond to flat vector bundles.

\begin{theorem}\label{fund-group-char-0}
Let $X$ be a smooth complex projective variety. Let us assume that $X$ is { D-quasi-affine}. 
Then  $\pi_1^{\et}(X)=0$. Moreover, all left $\cD_X$-modules, which are 
coherent as $\cO_X$-modules, are direct sums of finitely many $\cD_X$-modules 
isomorphic to the canonical $\cD_X$-module $\cO_X$.
\end{theorem}

\begin{proof}
Since $\pi_1^{\et}(X)$ is a profinite group, if $\pi_1^{\et}(X)\ne 0$ then there exists a non-trivial finite group $G$ and a surjective morphism $\pi_1^{\et}(X)\to G$. Taking, e.g., a regular representation $k[G]$ of $G$ we get a non-trivial linear representation of $\pi_1^{\et}(X)$. Since any representation of a finite group is unitary and it splits into a direct sum of irreducible representations, there exists also a non-trivial irreducible unitary representation $\pi_1^{\topo}(X)\to \GL (V)$ in some complex vector space $V$. 
The Riemann--Hilbert correspondence associates to this representation a vector bundle $E$ 
with an integrable connection $\nabla$. Since the representation is unitary, the stable Higgs bundle corresponding to $(E, \nabla)$ via Simpson's correspondence is simply $E$
with the zero Higgs field. In particular, $E$ is slope stable (with respect to any ample polarization) as a torsion free sheaf.  Since $\nabla$ corresponds to a left $\cD_X$-module structure on $E$, D-affinity of $X$ implies that 
$\Gamma (X, E)\ne 0$. But we know that $E$ has vanishing rational Chern classes (since $E$ carries a flat connection), so any non-zero section gives a map $\cO_X\to E$, which must be an isomorphism as $E$ is stable of degree $0$. But $\cO_X$ carries only one connection, as $h^0(X, \Omega_X)=h^1(X, \cO_X)=0$.
So $E$ is isomorphic to the $\cD_X$-module corresponding to $(\cO_X,d)$.
Then the corresponding representation is trivial, a contradiction. This shows that  $\pi_1^{\et}(X)=0$. 

A well-known result due to Malcev \cite{Ma} and Grothendieck \cite{Gr} shows that there are no nontrivial flat bundles on $X$. More precisely, since the topological fundamental group $\pi_1^{\topo}(X)$ is finitely generated and its profinite completion $\pi_1^{\et}(X)$ is trivial, by \cite[Theorem 1.2]{Gr} all finite dimensional  representations of $\pi_1^{\topo}(X)$ are also trivial. But by the Riemann--Hilbert correspondence such representations correspond to flat vector bundles,
so all flat vector bundles  are trivial (i.e., isomorphic to a direct sum of factors isomorphic to $(\cO_X,d)$).  
This shows the second part of the theorem.
\end{proof}

\begin{remark}
Proof of the first part of Theorem \ref{fund-group-char-0} can be obtained also in another way that we sketch here.
Namely, if $f: Y\to X$ is a finite \'etale covering and $Y$ is connected then $f_*\cO_Y$ is numerically flat. Therefore $E=f_*\cO_Y/ \cO_X$ is also numerically flat, so it admits a $\cD_X$-module structure. 
{  The short exact sequence 
$$0\to \cO_X\to f_*\cO_Y\to E\to 0$$
gives a short exact sequence 
$$0\to \Gamma (X,\cO_X)\to \Gamma (X, f_*\cO_Y)\to \Gamma (X, E)\to 0.$$
In characteristic zero this follows from the fact that the map $\cO_X\to f_*\cO_Y$ is split. 
Hence we have $\Gamma (X, E)=0$.  If $X$ is D-quasi-affine this} shows that $E=0$, so $f$ is a trivial covering. A similar argument works also in the positive characteristic case except that we
need $H^1(\cO_X)=0$ to assure that the sequence of sections is exact. We decided to give different arguments in both cases for two reasons. The first one is that the above argument in characteristic zero seems to give more insight into the proof (cf. proof of Proposition \ref{uniruled}). The second reason is that in positive characteristic this argument gives Theorem \ref{fund-group-char-p} only if one uses difficult  
\cite[Theorem 1.1]{EM}. In our proof of Theorem \ref{fund-group-char-p} we do not need to use this result.
\end{remark}

{ 
\begin{remark}
In Theorem \ref{fund-group-char-0} we assume only D-quasi-affinity of $X$. If one assumes that $X$ is D-affine as in
Theorem \ref{simply-connected},  then the proof of vanishing of  $\pi_1^{\et}(X)$ can be somewhat simplified (cf. proof of Corollary \ref{split-simple-connected}).
 \end{remark}
}

\subsection{Theorem \ref{simply-connected} in the positive characteristic case}

In this subsection we prove parts 1 and 2 of Theorem \ref{simply-connected} in case the base field $k$ has positive characteristic. As before, { 3 follows} from 1 and Proposition \ref{no-maps-to-curves}. 

\medskip

Let $X$ be a smooth variety defined over an algebraically closed field of positive characteristic. 
A \emph{stratified bundle} $\{E_n, \sigma_n\} _{n\in \ZZ_{\ge 0}}$ on $X$ is a sequence of locally free $\cO_X$-modules $E_n$ of finite rank and $\cO_X$-isomorphisms $\sigma_n: F_X^* E_{n+1}\stackrel{\simeq}{\to} E_n$.  Let us recall that by Katz's theorem \cite[Theorem 1.3]{Gi} the category of $\cD_X$-modules 
that are coherent as $\cO_X$-modules is equivalent to the category of stratified bundles.

\begin{theorem}\label{fund-group-char-p}
Let $X$ be a smooth projective variety defined over an algebraically closed field of positive characteristic. Let us assume that $X$ is { D-quasi-affine and $H^1(X, \cO_X)=0$}. Then  $\pi_1^{\et}(X)=0$. Moreover, all left $\cD_X$-modules, which are coherent as $\cO_X$-modules, are direct sums of finitely many $\cD_X$-modules 
isomorphic to the canonical $\cD_X$-module $\cO_X$.
\end{theorem}

\begin{proof} 
It is sufficient to prove that every stratified  bundle $E=\{E_n, \sigma_n\} _{n\in \ZZ_{\ge 0}}$
is a direct sum of the stratified bundles isomorphic to the stratified bundle $\bone_X$, corresponding to the $\cD_X$-module  $\cO_X$.

$D$-affinity of $X$ implies that $\Gamma (X, E_0)\ne 0$. Let us fix some integer $m\ge 0.$ Since $E(m):=\{E_{n+m}, \sigma_{n+m}\} _{n\in \ZZ_{\ge 0}}$ is a stratified bundle, we also have $\Gamma (X, E_m)\ne 0$ for all $m\ge 0$. By \cite[Proposition 2.3]{EM} there exists some $m\ge 0$ such that $E(m)$ is a successive extension of stratified bundles $U=\{U_n, \tau _n\} _{n\in \ZZ_{\ge 0}}$ such that all $U_n $ are slope stable of slope zero. By the same arguments as above  $\Gamma (X, U_n)\ne 0$, so $U_n\simeq \cO_X$. 
But  the sequence $\{\cO_X\}_{n\ge 0}$ admits only one structure of a stratified bundle (up to an isomorphism of stratified bundles), so $U\simeq \bone _X$.
Since $H^1(\cO_X)=0$ \cite[proof of Theorem 15]{dS} shows that   $E(m)$ is  a direct sum of  stratified bundles
isomorphic to $\bone _X$. But then $E$ is also a direct sum of  stratified bundles isomorphic to $\bone _X$.
This proves the second part of the theorem. Now equality $\pi_1^{\et}(X)= 0$ follows from  \cite[Proposition 13]{dS}.
\end{proof}
\medskip

\subsection{Maximal rationally connected fibrations of D-affine varieties}\label{MRC}

In this subsection we study maximal rationally connected fibrations of D-affine varieties in the characteristic zero case. First we prove a generalization of Proposition \ref{uniruled} that allows us to deal with rational maps.

Let $X$ be a smooth complete variety defined over an algebraically closed field $k$ of characteristic $0$.
Let $Y$ be a normal projective variety defined over $k$ and let $X^0\subset X$ and $Y^0\subset Y$ be non-empty open subsets. Let $f: X^0\to Y^0$ be a morphism such that $f_*\cO_{X^0}=\cO_{Y^0}$ (we do not require $f$ to be proper).

\begin{proposition}\label{rational-fibrations}
If $X$ is D-affine then one of the following holds:
\begin{enumerate}
\item $Y$ is uniruled, or
\item $Y\backslash Y^0$ has codimension  $1$ in $Y$.
\end{enumerate}
\end{proposition}

\begin{proof}
Since $Y$ is normal, it has singularities in codimension $\ge 2$ and hence without loss of generality we can assume that $Y^0$ is smooth, shrinking it if necessary. Let $i: X^0\hookrightarrow X$ and $j: Y^0\hookrightarrow Y$ denote the open embeddings. 

The proof is similar to that of Proposition \ref{uniruled}. Namely, let us fix an ample line bundle $L$ on $Y$ and consider $M=\cD_{Y^0}\otimes _{\cO_{Y^0}} j^*L^{-1}$. Then  $f^*M$ admits a left $\cD_{X^0}$-module structure. By \cite[Example 1.5.22 and Proposition 1.5.29]{HTT}  $i_*f^*M$ admits a left $\cD_{X}$-module structure. As in the proof of Proposition \ref{uniruled}, $\Gamma (X^0, f^*M) =\Gamma (X, i_*f^*M)\ne 0$ implies that $\Gamma (Y^0, S^m T_{Y^0}\otimes  _{\cO_{Y^0}} j^*L^{-1})\ne 0$ for some positive integer $m$.

If  $Y\backslash Y^0$ has codimension $\ge 2$ in $Y$ then by Hironaka's strong resolution of singularities
there exists a projective birational morphism $\pi: \tilde Y\to Y$ such that $\tilde Y$ is smooth,
$E=\pi^{-1}(Y\backslash Y^0)$ has pure codimension $1$ and $\pi $ is an isomorphism outside of $Y^0$. 
Let $\tilde j: Y^0 \hookrightarrow \tilde Y$ denote the lifting of $j$.
By construction we have $\Gamma (\tilde Y, \tilde j_* (S^m T_{Y^0})\otimes  _{\cO_{\tilde Y}} \pi^*L^{-1})\ne 0$ for some positive integer $m$. This implies that there exists some non-negative integer $n$ such that
$\tilde L:=\pi^*L(-nE)$ is a subsheaf of $S^mT_{\tilde Y}$. Let us note that
$$\tilde L(\pi^*L)^{d-1}=(\pi^*L)^{d}=L^d>0,$$
where  $d=\dim Y$. Let $A$ be an ample line budle on $\tilde Y$. Then for small $\epsilon>0$ we also have
$\tilde L(\pi^*L+\epsilon A)^{d-1}>0$. Hence $\mu_{\max, H} (S^mT_{\tilde Y})>0$, where 
$H=\pi^*L+\epsilon A$ is an ample divisor. As in the proof of Proposition \ref{uniruled} this implies that $\tilde Y$ is uniruled. Hence $Y$ is also uniruled.
\end{proof}

\medskip

Let $X$ be a smooth complete variety defined over an algebraically closed field $k$ of characteristic $0$.
Let $f: X\dasharrow Y$ be the maximal rationally connected fibration (see 
\cite[Chapter IV, Theorem 5.4]{Ko}). By definition there exist open subsets $X^0\subset X$ and $Y^0\subset Y$
and a morphism $f: X^0\to Y^0$ such that $f_*\cO_{X^0}=\cO_{Y^0}$ (note that we do not require $f: X^0\to Y^0$ to be proper). We assume that $Y$ is normal and projective (this can be always achieved by passing, if necessary, to another birational model of $Y$ using Chow's lemma and taking normalization).

\begin{proposition}\label{rat-connected}
If $X$ is D-affine then one of the following holds:
\begin{enumerate}
\item $X$ is rationally connected, or
\item $Y\backslash Y^0$ has codimension 
$1$ in $Y$,  $\dim X>\dim Y\ge 2$ and $\pi^1_{\et} (Y)=0$.
\end{enumerate}
\end{proposition}

\begin{proof}
If $\dim Y=0$ then $X$ is rationally connected, so we can assume that $\dim Y>0$. Then
the Graber--Harris--Starr theorem (see \cite[Corollary 1.4]{GHS}) implies that $Y$ is not uniruled.
Hence  by Proposition \ref{rational-fibrations} the complement of $Y^0$ in $Y$ has codimension 
$1$.  Let us note that $\dim X>\dim Y$ because $X$ is uniruled by Proposition \ref{uniruled}. If $\dim Y=1$ then $Y=\PP^1$ by Theorem \ref{fund-group-char-0} and Proposition \ref{no-maps-to-curves}. But this contradicts the Graber--Harris--Starr theorem, so $\dim Y\ge 2$. 

To see the last part we proceed as in the proof of Proposition \ref{no-maps-to-curves}. Namely, 
we can find a normal projective variety $\tilde X$, a birational morphism $\tilde X\to X$, which is an isomorphism over $X^0$, and a morphism $\tilde f: \tilde X\to Y$. As before $\tilde X$ is algebraically simply connected. Let us consider the Stein factorization of $\tilde f$
$$\tilde f: \tilde X\mathop{\longrightarrow}^{g} \tilde Y\mathop{\longrightarrow}^{h}  Y.$$
Note that by definition of a maximal rationally connected fibration, $f$ is proper over an open subset of $Y$ and then $h$ is an isomorphism over this subset. But $h$ is a finite birational morphism 
over a normal variety and hence it is an isomorphism. It follows that $\tilde f$ is a fibration.
Then we have a surjective map $\pi_1^{\et}(\tilde X)\to \pi_1^{\et}(Y)$ and hence $\pi_1^{\et} (Y)=0$.
\end{proof}

\begin{remark}
The above proposition strongly suggests that smooth projective D-affine varieties in characteristic zero are 
rationally connected. { In the case of  $3$-folds this fact follows from Proposition \ref{3-folds}.}
\end{remark}

\section{Proof of Theorem \ref{main-fibrations}}

We start with the following lemma that proves the first part of Theorem \ref{main-fibrations}.

\begin{lemma} \label{contractions}
Let $X$ be a smooth complete variety defined over an algebraically closed field $k$
and let $f:X\to Y$ be a fibration. Assume that $X$ is D-affine. Then for any effective divisor $D$
on $X$ the codimension of $f(D)$ in $Y$ is at most $1$. In particular, 
$Y$ does not admit any divisorial contractions.
\end{lemma}

\begin{proof}
Let us set $V=Y\backslash f(D)$ and $U=f^{-1} (V)$. If $f(D)$ has codimension $\ge 2$ in $Y$ then
$$k=\Gamma (Y, \cO_Y)\mathop{\to}^{\simeq} \Gamma (V, \cO_V)\mathop{\to}^{\simeq}  \Gamma (U, \cO_U).$$
Since $D\subset X\backslash U$ this contradicts Lemma \ref{Thomsen}.

Now assume that $Y$ admits a birational morphism  $g: Y\to Z$ onto a normal variety $Z$ and 
the exceptional locus $E$ of $g$ has codimension $1$. Then $D=f^{-1}(E)$ contains a divisor and 
$h=gf: X\to Z$ is a fibration. But $h(D)$ has codimension $\ge 2$ in $Z$, a contradiction.
\end{proof}

\begin{remark}
\begin{enumerate}
\item If $g: Y\to Z$ is a birational morphism, $Z$ is normal and locally $\QQ$-factorial then the exceptional locus has pure codimension $1$, so we can apply the above corollary.
\item Let $X$ be as in Lemma \ref{contractions} and let $f:X\to Y$ be a fibration over a smooth variety $Y$. If $X$ is D-affine then by \cite[Corollary 6.12]{Ar} the above lemma implies that $Y$ does not contain any smooth divisors with ample conormal bundle.
\end{enumerate}
\end{remark}

The following proposition proves part 2 of Theorem \ref{main-fibrations}. It is a small generalization of \cite[2.14 Proposition (a)]{HP} that follows Kashiwara's  proof of \cite[Theorem 1.4.1]{Ka} (which is based on an idea used by Beilinson-Bernstein in proof of their theorem). However, the result is stated in characteristic zero and in positive characteristic it needs to be reformulated. Even in characteristic zero checking (in the notation of \cite{HP}) that $\cD^f=\cD_X$  for $\cD=\cD_Y$  requires a non-trivial computation that is missing in \cite{HP}. Since the authors only sketch the arguments and add some unnecessary assumptions, we give a full proof of the result.

\begin{proposition}\label{vanishing-for-fibers}
Let $X$ a smooth complete variety defined over an algebraically closed field $k$.
Assume that $X$ is almost D-affine. Let $f: X\to Y$ be a surjective morphism onto a projective variety $Y$.
Then for any $\cD_X$-module $M$ we have $R^jf_*M=0$ for all $j>0$. In particular, we have $R^jf_*\cO_X=0$ for all $j>0$ and if $F$ is a general fiber of $f$ then $H^j(F, \cO_F)=0$ for all $j>0$.
\end{proposition}

\begin{proof}
Let $L$ be any globally generated line bundle on $X$. Then the surjection $\cO_X\otimes _k \Gamma (X, L)\to L$
induces a surjection $L^{-1}\otimes _k \Gamma (X, L)\to \cO_X$. After tensoring with $\cD_X$ we get
a surjective map of (left) $\cD_X$-modules
$$\varphi: \cD_X\otimes _{\cO_X} (L^{-1}\otimes _k \Gamma (X, L))\to \cD_X.$$
Let us note that if $X$ is D-affine then this map has a section (as a map of $\cD_X$-modules).
This follows from the fact that
$$\Hom _{\cD_X}(\cD_X,  \cD_X\otimes _{\cO_X} (L^{-1}\otimes _k \Gamma (X, L)))=\Gamma (X,  \cD_X\otimes _{\cO_X} L^{-1})\otimes _k \Gamma (X, L)\to \Hom _{\cD_X}(\cD_X,\cD_X)=\Gamma (X,  \cD_X)$$
 is surjective, as its cokernel is contained in $H^1(X, \ker \varphi)=0$. 
 
 Using $\cHom _{\cD_X} ( \bullet, \cD_X)$ we get a split map of right $\cD_X$-modules
$$\cD_X\to (L\otimes _k \Gamma (X, L)^*) \otimes _{\cO_X} \cD_X.$$
Taking $\bullet \otimes _{\cD_X} M$  we get  a split map 
$$M\to (L\otimes _k \Gamma (X, L)^*) \otimes _{\cO_X} M$$
of sheaves of abelian groups.

 Now let us take an ample line bundle $A$ on $Y$. Let us consider a coherent $\cO_X$-submodule $G$ of $M$.  
 Then for large $m\gg 0$ $L:=f^* A^{\otimes m}$ is globally generated  and $R^jf_*(G\otimes L)=0$ for all $j>0$.
 We have a commutative diagram
 $$
\xymatrix{
R^jf_*G \ar[d]\ar[r]&R^jf_*( (L\otimes _k \Gamma (X, L)^*) \otimes _{\cO_X} G) =0\ar[d]\\
R^jf_*M\ar[r]& R^jf_*( (L\otimes _k \Gamma (X, L)^*) \otimes _{\cO_X} M) \\
}
 $$
in which the lower horizontal map is split. So the map $R^jf_*G\to R^jf_*M$ is zero.
Since $R^jf_*M$ is the direct limit of $R^jf_*G$, where $G$ ranges over all coherent  $\cO_X$-submodules $G$ of $M$, we get  required vanishing of $R^jf_*M$. Applying this to $M=\cO_X$ we get the last part of the proposition.
\end{proof}

\begin{remark} \label{vanishing-H1}
If under the assumptions of Proposition \ref{vanishing-for-fibers} the morphism
$f: X\to Y$ is a fibration (i.e., $f_*\cO_X=\cO_Y$) then the Leray spectral sequence implies 
that $H^i(Y, \cO_Y)=0$ for $i>0$.
\end{remark}

By  $D _{qc}(\Mod{\cO_X})$  we denote the full  subcategory of the (unbounded) derived category
of the category of $\cO_X$-modules, consisting of complexes whose cohomology sheaves are quasi-coherent. 

The above proposition implies the following corollary:

\begin{corollary}
In the notation of Proposition \ref{vanishing-for-fibers} for any $\cD_X$-module $M$ 
 the canonical map $f_*M\to Rf_*M$ is an isomorphism in  $D _{qc}(\Mod{\cO_Y})$. In particular, the
functor $f_*: \Mod{\cD_X}\to \Mod{\cO_Y}$ is exact.
\end{corollary}
\medskip

In the above corollary $f_*$ denotes the composition of the forgetful functor $\Mod{\cD_X}\to \Mod{\cO_X}$
with  the direct image $f_*: \Mod{\cO_X}\to \Mod{\cO_Y}$. 

\medskip

Let us recall that a sheaf $G$ on a scheme of pure dimension is called \emph{Cohen--Macaulay} if for every point $x\in X$, the depth of $G$ at $x$ is equal to the codimension of $x$ in $X$ (see \cite[Definition 11.3]{Kol2}).
By \cite[Proposition 3.12]{Kol} Proposition \ref{vanishing-for-fibers} implies the following result:

\begin{corollary}\label{Kollar}
Let $X$ be a smooth projective variety defined over an algebraically closed field $k$ of characteristic zero.
Assume that $X$ is almost D-affine. If $f: X\to Y$ is any surjective morphism onto some normal projective variety $Y$ and $f_*\cO_X$ is torsion free (e.g., $f$ is a fibration) then $Y$ has only rational singularities and $f_*\cO_X$ is
a Cohen--Macaulay sheaf.
\end{corollary}

\medskip

\begin{example}\label{Lauritzen}
To show further usefulness of Proposition \ref{vanishing-for-fibers}  let us reprove Lauritzen's result that some unseparated flag varieties are not D-affine (see \cite[Section 4]{Lau}). Namely, let $n\ge 2$ and let $X$ be the zero scheme 
$x_0y_0^{m}+...+x_ny_n^{m}=0$ in $\PP^n\times \PP^n$.  Let $f: X\to \PP ^n$ be the projection onto the first factor and let $F$ be any fiber of $f$.
A short exact sequence
$$0\to \cO_{\PP^n}(-m)\to \cO_{\PP^n}\to \cO_F\to 0$$ 
shows that if $m\ge n+1$ then $H^{n-1}(F,\cO_F)\simeq H^{n}(\PP^r, \cO_{\PP^n}(-m))\ne 0$, so $X$ is not D-affine.

If $X$ is considered over an algebraically closed field of characteristic $p$ and $m=p^r$ for some $r\ge 1$
then $X$ is an unseparated flag variety. In this case all fibers $F$ of $f$ are multiplicity $p^r$ hyperplanes in $\PP ^n$. By the above $X$ is not D-affine if $p^r\ge n+1$.
\end{example}

\medskip

{

The following proposition proves part 3 of Theorem \ref{simply-connected}. It shows that any surjective morphism from D-affine variety behaves like a flat morphism (for $D$-modules).

\begin{proposition}\label{technical}
Let $X$ a smooth complete D-affine variety defined over an algebraically closed field $k$.
Let $f: X\to Y$ be a surjective morphism onto a smooth projective variety $Y$.
Let us assume that $f_*\cO_X$ is locally free. Then for any $\cD_Y$-module $N$ the 
following conditions are satisfied:
\begin{enumerate}
\item the canonical map $Lf^*N\to f^*N$ is an isomorphism in  $D _{qc}(\Mod{\cO_X})$,
\item $H^j(Y,f_*\cO_X\otimes _{\cO_Y} N)=0$ for all $j>0$.
\end{enumerate}
In particular, the functor $f^*: \Mod{\cD_Y}\to \Mod{\cD_X}$ is exact.
\end{proposition}

\begin{proof}
Let us recall that for any bounded complex $C^{\bullet}$ of quasi-coherent $\cO_X$-modules we have a spectral sequence 
$$E^{ij}_2=H^i(X, \cH ^j (C^{\bullet}))\Rightarrow H^{i+j}(C^{\bullet}).$$
Let $N$ be any $\cD_Y$-module. Let us recall that by assumption $N$ is quasi-coherent as an 
$\cO_Y$-module. Applying the above spectral sequence to $Lf^*N$ (which is represented by a bounded complex of quasi-coherent $\cO_X$-modules) we get 
$$E^{ij}_2=H^i(X, L^jf^*N)\Rightarrow H^{i+j}(X, Lf^*N).$$
Note that $L^jf^*N$ are quasi-coherent $\cO_X$-modules carrying a left $\cD_X$-module structure
(see \cite[1.5]{HTT} for the characteristic $0$ and \cite[Section 2]{Ha2} for the positive characteristic case). Hence by D-affinity of $X$ we have $H^i(X, L^jf^*N)=0$ for all $i>0$ and any $j$. 
So the above spectral sequence degenerates to $$H^0(X, L^jf^*N)=H^j(X, Lf^*N)$$
for all $j$.

By Proposition \ref{vanishing-for-fibers} the canonical map $f_*\cO_X\to Rf_*\cO_X$ is an isomorphism
in  $D _{qc}(\Mod{\cO_Y})$. Hence by the projection formula  (see \cite[Proposition 5.3]{Ne}) there exist natural isomorphisms 
$$f_*\cO_X\otimes _{\cO_Y}^L N\mathop{\to}^{\simeq} Rf_*\cO_X\otimes ^L_{\cO_Y} N\mathop{\to}^{\simeq} Rf_*(Lf^*N)$$
in  $D _{qc}(\Mod{\cO_Y})$.

Since $f_*\cO_X$ is a locally free $\cO_Y$-module, the canonical map  $f_*\cO_X\otimes _{\cO_Y}^L N\to f_*\cO_X\otimes _{\cO_Y} N$ is an isomorphism. 
 In particular, $R^0f_*(Lf^*N)=f_*\cO_X\otimes _{\cO_Y} N$ and
$R^{j}f_*(Lf^*N)=0$ for $j\ne 0$. 
Hence the Leray spectral sequence 
$$E^{ij}_2=H^i(Y,R^jf_*(Lf^*N))\Rightarrow H^{i+j}(X,Lf^*N)$$
degenerates to 
$$H^i (Y, f_*\cO_X\otimes _{\cO_Y} N)=H^i (Y, R^0f_*(Lf^*N))=H^i (X, Lf^*N).$$
It follows that 
$$H^0(X, L^jf^*N)=H^j(X, Lf^*N)=H^j (Y, f_*\cO_X\otimes _{\cO_Y} N)$$
for all $j$. But then $H^0(X, L^jf^*N)=0$ for $j<0$ and D-affinity of $X$ implies that 
$L^{j}f^*N=0$ for all $j<0$. Since $L^{j}f^*N=0$ for all $j< 0$, we have a natural isomorphism $Lf^*N\stackrel {\simeq}{\to} f^*N$ in  $D _{qc}(\Mod{\cO_X})$.
We also get $H^j (Y, f_*\cO_X\otimes _{\cO_Y}N)=0$ for all $j>0$ as $L^jf^*N=0$ for $j>0$.
\end{proof}

\medskip

The following lemma explains the meaning of the assumption in Proposition \ref{technical}. 

\begin{lemma}
Let $f: X\to Y$ be a projective morphism of smooth $k$-varieties and let 
$X\stackrel{f'}{\to}Y'\stackrel{g}{\to}Y$ be Stein's factorization of $f$.
Then  $f_*\cO_X$ is locally free if and only if $Y'$ is Cohen--Macaulay.
\end{lemma}

\begin{proof}
$Y'$ is Cohen--Macaulay if and only if $\cO_{Y'}$ is  Cohen--Macaulay.
 By definition of Stein's factorization we have $g_*\cO_{Y'}=f_*\cO_X$. 
 Hence by \cite[Lemma 11.4]{Kol2} $Y'$ is Cohen--Macaulay  if and only if
$f_*\cO_X$ is  Cohen--Macaulay, which by the same lemma is equivalent to 
 $f_*\cO_X$ being locally free.
\end{proof}

\medskip

\begin{remark}
Since Cohen--Macaulay sheaves on smooth varieties are locally free (see \cite[Lemma 11.4]{Kol2}),
Corollary \ref{Kollar} implies that if $k$ has characteristic zero and $f_*\cO_X$ is torsion free then  $f_*\cO_X$ is locally free.
\end{remark}

The first part of the above proposition implies the following corollary:

\begin{corollary}
In the notation of Proposition \ref{technical} the functor $f^*: \Mod{\cD_Y}\to \Mod{\cD_X}$ is exact.
\end{corollary}

\begin{corollary} \label{split-simple-connected}
If, in the notation of Proposition \ref{technical}, the canonical map $\cO_Y\to f_*\cO_X$ is split in $\Mod{\cO_Y}$  then  $Y$ is almost D-affine. In particular, $Y$ is algebraically simply connected.
\end{corollary}

\begin{proof}
The first assertion is clear as for any $\cD_Y$-module $M$, $H^i(Y, M)$ is a direct summand of $H^i(Y,f_*\cO_X\otimes _{\cO_Y} M)$.
To prove the second assertion note that if $M$ is a $\cD_Y$-module, which is coherent as an $\cO_Y$-module,
then it has vanishing numerical Chern classes. Moreover, such $M$ is locally free of finite rank $r$ 
(as an $\cO_Y$-module). Since  $H^i(Y, M)=0$ and  $H^i(Y, \cO_Y)=0$ for any $i>0$, 
the Riemann--Roch theorem implies that 
$$\dim \, \Gamma (Y, M)=\chi (Y, M)=r\chi (Y, \cO_Y)=r.$$
It follows that $M\simeq \cO_Y^r$. This implies vanishing $\pi_1^{\et} (Y)=0$.
\end{proof}

\begin{remark}
The above corollary implies that if $f: X\to Y$ is an \'etale morphism of smooth projective varieties and $X$
is D-affine then $f$ is an isomorphism. This is related to generalized Lazarsfeld's problem \cite[p.~59]{Laz}
asking if smooth images of flag varieties under finite morphisms are flag varieties.
\end{remark}

\begin{corollary} \label{image-of-D-affine-is-D-affine}
In the notation of Proposition \ref{technical} let us assume that $f_*\cO_X=\cO_Y$ (i.e., $f$ is a fibration).
Then $Y$ is D-affine.
\end{corollary}

\begin{proof}
By Proposition \ref{technical}  for any $D_Y$-module $M$ we have $H^j (Y, M)=0$ for all $j>0$.
By  \cite[Proposition 1.5.2]{Ka} (or \cite[1.4]{HTT}) to finish the proof of D-affinity of $Y$, it is sufficient to show that if $M\ne 0$ then $\Gamma (Y, M)\ne 0$.  By Proposition \ref{vanishing-for-fibers} the canonical map $\cO_Y\to Rf_*\cO_X$ is an isomorphism in  $D _{qc}(\Mod{\cO_Y})$.
Hence  by Proposition \ref{technical} and the projection formula there exist natural isomorphisms 
$$M\mathop{\to}^{\simeq} Rf_*\cO_X\otimes ^L_{\cO_Y} M\mathop{\to}^{\simeq} Rf_*(Lf^*M)\mathop{\to}^{\simeq} Rf_*(f^*M)$$
in $D _{qc}(\Mod{\cO_Y})$. So the canonical map  $M\to f_*(f^*M)$ of $\cO_Y$-modules is an isomorphism and $R^if_*(f^*M)=0$ for all $i>0$. In particular, if $M\ne 0$ then $f^*M\ne 0$. 
Hence by D-affinity of $X$ we have
$$\Gamma (Y, M)=\Gamma (Y, f_*f^*M)=\Gamma (X, f^*M)\ne 0,$$
which finishes the proof of D-affinity of $Y$. 
\end{proof}

\medskip

\begin{remark}
In his PhD thesis B. Haastert showed that if $f: X\to Y$ is a locally trivial fibration of smooth varieties with smooth D-affine fibers and $X$ is D-affine then $Y$ is D-affine (see \cite[Satz 3.8.9]{Ha-PhD}). The above proposition is a generalization of this fact without any assumptions on the fibers and local freeness of the fibration. Our proof is completely different.
\end{remark}

\begin{proposition}\label{BB-glob-gen}
Let $X$ a smooth complete variety defined over an algebraically closed field $k$.
Assume that $X$ is D-quasi-affine. Let $f: X\to Y$ be a surjective morphism onto a projective variety $Y$.
Then for any $\cD_X$-module $M$ the canonical map $\cD_X\otimes_{\cO_X} f^*f_*M\to M$ is surjective.
\end{proposition}

\begin{proof}
Let $g: Y\to \Spec k$ be the structure morphism and let $h=g f$. We have canonical maps
$$h^*(h_*M)\to f^*(g^*(g_*(f_*M)))\to f^*f_*M\to M.$$
Hence the surjection $\cD_X\otimes_k\Gamma(X, M) =\cD_X\otimes_{\cO_X} h^*h_*M\to M$ factors through 
$\cD_X\otimes_{\cO_X} f^*f_*M\to M$. 
\end{proof}

}

\section{Uniruledness of surfaces in positive characteristic}

Let us recall the following uniruledness criterion of Miyaoka and Mori 
(see \cite[Corollary 3]{MM}).

\begin{theorem}\label{Miyaoka-Mori}
Let $X$ be a (possibly non-normal) $\QQ$-Gorenstein projective variety of dimension $n$
defined over an algebraically closed field of any characteristic. If there exist ample 
divisors $H_1,...,H_{n-1}$ such that $K_XH_1...H_{n-1}<0$ then $X$ is uniruled.
\end{theorem}

Miyaoka used this criterion to prove that if $X$ is a smooth projective variety defined over an algebraically closed field of characteristic zero then the cotangent bundle of $X$ is generically semipositive unless $X$ is uniruled (see \cite[Corollary 8.6]{Miy}).
Unfortunately, Miyaoka's criterion of uniruledness does not work in positive characteristic (see the next subsection). However,  we give a certain generalization of this criterion that works 
in the surface case in an arbitrary characteristic (see Theorem \ref{symmetric-power}).

\subsection{Ekedahl's example revisited} \label{Ekedahl-example}

Here we give an example of a non-uniruled surface of general type $\tilde X$ such that 
$\Omega_{\tilde X}$ is not generically semipositive for some ample polarizations. 
In particular, such $\tilde X$ has an unstable tangent bundle.

\medskip

Let $Y$ be a smooth projective surface defined over a field $k$ of characteristic $p$.
Let $L$ be a very ample line bundle on $Y$ and let $s\in H^0(Y, L^p)$ be a general section.
Let $\pi: \VV (L)\to Y$ be the total space of $L$, i.e., $\VV (L)=\Spec _Y \bigoplus _{i\ge 0} L^{-i}$.
Since by the projection formula $\pi_*(\pi^*L)=\bigoplus _{i\ge -1} L^{-i}$, $\pi^*L$ has a canonical section $t_L$
corresponding to $1\in H^0(X, \cO_X)$. Therefore both $t_L^p$ and $s$ can be treated as sections of $\pi^*L^p$.

Let $\varphi: X=Y[\sqrt[p]{s}]\to Y$ be a degree $p$ cyclic covering defined by $s$, i.e., $X$ is defined as the zero set of $t_L^p-s$. Then we have an exact sequence
$$\pi^*L^{-p}\simeq \cO_{\VV(L)}(-X)\mathop{\to}^{d_X} \Omega^1_{\VV(L)}\to \Omega^1_X\to 0,$$
which on $X$ induces  an exact sequence
$$0\to  \varphi^*L^{-p}\to \varphi^*\Omega _Y\to \Omega_X \to \varphi ^*L^{-1}\to 0.$$
In fact, the first map is the pull-back of an $\cO_Y$-linear map $d s:L^{-p}\to \Omega_Y $ defined locally by the property that if $\tau \in L^{-1}$ is a local generator then $(ds)(\tau^{p})=d(s\tau^{p})$.

Our assumptions imply that $X$ is integral and $\varphi$ is purely inseparable, but $X$ is usually singular. It is singular exactly over the set of points where  $ds\in H^0(Y, L^p\otimes \Omega_Y)$ vanishes 
(such points are called \emph{critical}), as at such points the cokernel of $L^{-p}\to \Omega_Y$ is not locally free.

We always have  $K_X=\varphi^*(K_Y+(p-1)L)$ but $X$ can be non-normal in general (this happens if and only if the set of critical points of $ds$ has codimension $1$ in $Y$).

Now assume that $Y$ is an abelian surface and $L$ is $2$-jet generated (i.e., at every point $x\in X$ the evaluation map $H^0(X, L)\to L\otimes \cO_{X,x}/m_x^3$ is surjective).
Note that this last condition can be always arranged (e.g., if $H$ is ample and globally generated on $Y$ then $4H$ is $2$-jet generated).

Then $K_X$ is ample but $X$ is always singular because the cokernel of $ds$ cannot be locally free (this can be seen by a simple computation of $c_2(\Omega_Y)$).
However, our assumptions on $L$ imply that a general section of $H^0(Y, L)$ has only a finite number of nondegenerate critical points (see \cite[Chapter V, Exercise 5.7]{Ko}).
Therefore $X$ is normal and the cokernel of $ds$ is torsion free. Let us write the cokernel of $ds$ as $I_ZL^p$ for some $0$-dimensional scheme $Z$.

Let $f:\tilde X\to X$ be a resolution of singularities.
Then the canonical map $f^*\Omega_X\to \Omega_{\tilde X}$ is generically an isomorphism, so we have an injective map
$$f^*\varphi^*(I_ZL^{p})/\tors \hookrightarrow \Omega _{\tilde X}.$$ 
This shows that $\Omega_{\tilde X}$ has a quotient that is rank $1$ torsion free sheaf $\tilde L$ with first Chern class of the form $f^*\varphi^*L^{-1}+\hbox{($f$-exceptional divisor)}$
and 
$$\mu_{f^*\varphi^*L}(\Omega_{\tilde X})=\frac{p-1}{2}L^2>0>\mu_{f^*\varphi^*L} (\tilde L)=\mu_{f^*\varphi^*L}(f^*\varphi ^*L^{-1})=-L^2.$$ 
Clearly, the same inequalities hold for ample polarizations of the form $f^*\varphi^*L -\hbox{(small $f$-exceptional divisor)}$. 
Therefore $\tilde X$ is a surface of general type and $\Omega_{\tilde X}$ is not generically $H$-semipositive
for some ample $H$, even though $\tilde X$ is non-uniruled (as it admits a generically finite map onto an abelian surface).

\begin{remark}\label{Eke-rem}
This example is a corrected version of Ekedahl's example as described by Miyaoka in  \cite[Example 8.8]{Miy} (see also \cite[p.~145--146]{Ek} for a similar example but with a different aim in mind).
The example in \cite{Miy} does not work as stated due to existence of singularities of the covering.
\end{remark}

\subsection{ Generic semipositivity of cotangent bundle in positive characteristic}

The following theorem is a positive characteristic { version } of Miyaoka's generic semipositivity theorem.

\begin{theorem} \label{symmetric-power}
Let $X$ be a smooth projective surface defined over an algebraically closed field 
$k$ of characteristic $p$. Let us fix an ample divisor $H$ on $X$.
{
If $X$ is not uniruled then either $\Omega_X$ is generically $H$-semipositive or
$T_X$ is not slope $H$-semistable and $0<\mu_{\max, H}(T_X)\le K_XH/(p-1)$.
}
\end{theorem}

\begin{proof} {
By Theorem \ref{Miyaoka-Mori} we have $K_XH\ge 0$. Let us assume that $\Omega_X$ is not generically $H$-semipositive. Then  $L_{\max, H}(T_X)=-L_{\min, H}(\Omega_X)> 0$ and there exists some $m\ge 0$ such that $\mu_{\max, H}((F_X^m)^*T_X)>0$.
Then $(F_X^m)^*T_X$ is not slope semistable as $\mu ((F_X^m)^*T_X)=-p^mK_XH\le 0$
and the maximal destabilizing subsheaf $L\hookrightarrow (F_X^m)^*T_X$ is a line bundle with $LH>0$.
But $(F_X^m)^*T_X$ is a subsheaf of $S^{mp}T_X$, so it is sufficient to prove that 
if for some $j\ge 1$ the $j$-th symmetric power $S^jT_X$ of the tangent bundle 
contains a line bundle $L$  such that $LH>0$ then $T_X$ is not slope $H$-semistable and $0<\mu_{\max, H}(T_X)\le K_XH/(p-1)$.
}

By assumption there exists some $j>0$ such that $S^j T_X\otimes L^{-1}$ has a non-zero section.
This section gives rise to an effective divisor $D\in |\cO_{\PP (T_X)}(j)-\pi^*L|$, where 
$\pi : \PP(T_X)\to X$ denotes the projective bundle. Let us write $D$ as a sum
$$D=a_1D_1+...+a_mD_m,$$
where $a_i>0$ and $D_i$ are irreducible and reduced divisors. Since $\Pic (\PP(T_X))$
is generated by $\cO_{\PP (T_X)}(1)$ and $\pi^* \Pic X$, for each $i=1,...,m$ we can  
find $n_i\ge 0$ and line bundles $L_i$ on $X$ such that $D_i\in |\cO_{\PP (T_X)}(n_i)-\pi^*L_i|$.
By the adjunction formula we have
$$K_{D_i}=(K_{\PP (T_X)}+D_i)|_{D_i}=(\cO_{\PP (T_X)}(n_i-2)-\pi^*L_i)|_{D_i}.$$
Using the Leray--Hirsch formula we get
$$K_{D_i}\pi^*H= (n_i-2)n_i (-K_XH)-(2n_i-2)L_iH.$$
By the Miyaoka--Mori theorem if $K_{D_i}\pi^*H<0$ then $D_i$ is uniruled (note that $\pi^*H$ is only nef,
but we can always find an ample divisor $H_i$ on $D_i$ such that $K_{D_i}H_i<0$). If $n_i>0$ then $X$ is also uniruled, { a contradiction.  So we have $K_{D_i}\pi^*H\ge 0$ for all $i$ such that $n_i>0$.} 
Thus if $n_i\ge 2$ then we get 
$$L_iH \le \frac{(n_i-2)n_i}{2n_i-2} (-K_XH)\le 0.$$
If $n_i=1$ then $D_i$ gives a section of $T_X\otimes L_i^{-1}$, so $L_iH\le \mu _{\max} (T_X)$.
If $n_i=0$ then let us take some $m$ such that $\cO_{\PP (T_X)}(1)+m\pi^*H$ is ample.
Then 
$$0\le D_i\pi^*H(\cO_{\PP (T_X)}(1)+m\pi^*H)=-L_iH,$$
so $L_iH\le 0$. Therefore
$$0<LH=\sum a_iL_iH\le \left( \sum_{\{i: n_i=1\}} a_i \right) \mu _{\max} (T_X)$$
and hence $\mu_{\max} (T_X)> 0$.

Since $K_XH\ge 0$, $T_X$ is not slope $H$-semistable and the maximal destabilizing subsheaf $M\subset T_X$ 
has rank $1$. By assumption $MH>0$. Let us note that $\Hom (F_X^*M, T_X/M)=0$ since 
$\mu (F_X^*M)=pMH>0> \mu (T_X/M)=(-K_X-M)H$. Therefore $M\subset T_X$ defines a $1$-foliation.
Let $f: X\to Y$ be the quotient by this $1$-foliation. There exists an ample $\QQ$-divisor $H'$
on $Y$ such that $f^*H'=H$. Then we have
$$K_YH'=K_XH -(p-1)MH.$$
If $K_XH<(p-1)MH$ then $K_YH'<0$, so $Y$ is uniruled. But then $X$ is also uniruled, { a contradiction.
This shows that $(p-1)MH\le K_XH$.}
\end{proof}

\begin{remark}
In the example from the previous subsection $\tilde X$ is not uniruled, {  $\Omega_{\tilde X}$ is not generically $H$-semipositive}  and  $T_{\tilde X}$ contains a line bundle $M$ such that $MH=K_{\tilde X}H/(p-1)$ 
for some nef and big divisor $H$. This shows that the above upper bound on $\mu_{\max, H}(T_X)$ is optimal.
\end{remark}

{ 
\section{D-affine varieties in low dimensions}\label{surfaces}

\subsection{Surfaces that are images of D-affine varieties} 
}

The main aim of this section is to prove the following slightly more precise version of Theorem \ref{D-affine-surfaces}.

\begin{theorem}\label{classification-2}
Let $X$ be a smooth projective variety defined over an algebraically closed field $k$
and let $f: X\to Y$ be a fibration over a smooth projective surface { $Y$}.
If $X$ is D-affine then { $f$ is flat, $Y$ is D-affine and } one of the following holds:
\begin{enumerate}
\item $Y=\PP^2$,
\item $Y=\PP^1\times \PP^1$,
\item $2\le \chr k\le 7$ and $Y$ is an algebraically simply connected surface of general type with $p_g(Y)=q(Y)=0$ and $1\le K_Y^2\le 9$. Moreover, $Y$ is not uniruled and $K_Y$ is ample.
\end{enumerate}
\end{theorem}
 
\begin{proof}
{ By Lemma \ref{contractions} all fibers of $f$ have dimension $\le \dim X-2$. Since the dimension of any fiber  is at least $\dim X-2$, $f$ is equidimensional. So $f$ is also flat and by Proposition 
\ref{image-of-D-affine-is-D-affine}  $Y$ is D-affine. }
By Lemma \ref{contractions} and Artin's (or Grauert's if $k=\CC$) criterion of contractibility (see \cite[Corollary 6.12]{Ar}) $Y$ does not contain any irreducible curves $C$ with $C^2<0$. In particular, $Y$ is minimal. By Theorem \ref{simply-connected} we also know that $Y$ 
does not admit any maps onto curves of genus $\ge 1$.
 
Therefore if the Kodaira dimension $\kappa(Y)=-\infty$ then $Y$ is a minimal rational surface. By 
Lemma \ref{contractions} $Y$ is not the Hirzebruch surface $F_n$ for $n\ge 2$ as $F_n$
contains   a curve with self-intersection $(-n)$. Hence $Y=\PP ^2$ or $Y=\PP^1\times \PP^1$. These surfaces are D-affine in an arbitrary characteristic.

If $\chr k=0$ then by Proposition \ref{uniruled} we know that $Y$ is uniruled and hence $\kappa(X)=-\infty$.
So we can assume that $\chr k=p>0$. 

If the Kodaira dimension $\kappa(Y)\ge 0$ then $K_Y$ is nef, so $c_1(Y)^2\ge 0$.  D-affinity of { $Y$ implies $h^1(\cO_Y)=h^2(\cO_Y)=0$}, so $p_g(Y)=h^0(K_Y)=h^2(\cO_Y)=0$ and $\chi (Y, \cO_Y)=1$. 

If $c_1^2(Y)=0$ then \cite[E.4, Theorem]{Re} implies that $\kappa (Y)=1$. 
In this case $Y$ admits an elliptic or quasi-elliptic fibration $g: Y\to \PP ^1$.
In particular, we have $R^1g_*\cO_X\ne 0$, {  which contradicts 
Proposition \ref{vanishing-for-fibers}.}

If $c_1^2(Y)>0$ then $Y$ is a minimal surface of general type. Since $Y$ does not contain any $(-2)$ curves,
the canonical divisor $K_Y$ is ample. Since $h^1(\cO_Y)=0$ we have $b_1(Y)=0$ and  
$c_2(Y)=2-2b_1(Y)+b_2(Y)\ge 3$. Therefore Noether's formula gives $K_Y^2=12 -c_2(Y)\le 9$. 

{
Let us recall that $K_Y$ is ample and let us assume that $Y$ is not uniruled.
As in proof of Proposition \ref{uniruled} there exists some $j>0$ such that { $\Gamma ^j T_Y\otimes \omega_Y^{-1}$} has a non-zero section. Therefore by Theorem \ref{L_max-tensor-product}
$$0< K_Y^2\le  \mu_{\max, K_Y}(\Gamma ^j T_Y) \le L_{\max, K_Y}(\Gamma ^j T_Y)= j\, L_{\max, K_Y}(T_Y).$$
Then by  Theorem \ref{symmetric-power}} $T_Y$ contains a saturated line bundle $M$ such that  $K_Y^2\ge (p-1)MK_Y>0$. In particular, we have $K_Y^2\ge p-1$. But $K_Y^2\le 9$, which implies $p\le 7$.

Now let us assume that $Y$ is uniruled. Let us consider the maximal rationally chain connected fibration $g: Y\dasharrow Z$ (see \cite[IV.5]{Ko}). Since $Y$ is uniruled, we have $\dim Z\le 1$. If $\dim Z=0$ then $Y$ is rationally chain connected. But since $\dim Y=2$, $Y$ is rational, a contradiction. So $Z$ is a curve.
By Proposition \ref{no-maps-to-curves} $Z=\PP ^1$, which again implies that $Y$ is rational (by the analogue of L\"uroth's theorem in dimension $2$), a contradiction.
\end{proof}

\begin{remark}
Classification of algebraically simply connected surfaces of general type with $p_g(Y)=q(Y)=0$ is  a well-known open problem.  By Noether's formula we get $c_1^2(Y)+c_2(Y)=12$. Since $1\le c_1^2(Y)\le 9$ the possible number of families of such surfaces is very limited. The first examples of such surfaces were constructed by Barlow, but by construction the canonical divisor of such a surface is not ample, so Barlow's surfaces are not D-affine.

In \cite{LN} Lee and Nakayama showed that for any algebraically closed field $k$ and any $n\in \{1,2,3,4\}$, there exist algebraically simply connected minimal surfaces $Y$ of general type over $k$ with $K_Y^2=n$, $p_g(Y)=q(Y)=0$ and with  ample $K_Y$ (except possibly if $\chr k=2$ and $n=4$).  It is not clear how to check if these examples are D-affine in characteristics $2\le p\le 7$. The author does not know any other examples of smooth projective surfaces  of general type with $p_g(Y)=q(Y)=0$ that are algebraically simply connected.
\end{remark}

{

\subsection{Smooth projective D-affine $3$-folds}

The main aim of this section is to prove the following proposition:

\begin{proposition} \label{3-folds}
Let $X$ be a smooth projective variety of dimension $3$ defined over an algebraically closed field $k$.
If $X$ is D-affine then  one of the following holds:
\begin{enumerate}
\item $X$ is a smooth Fano $3$-fold with $b_2(X)=1$.
\item  There exists a fibration $f: X\to \PP^1$ such that every fiber of $f$ with reduced scheme structure is a del Pezzo surface.
\item There exists a smooth projective D-affine surface $Y$ and a flat conic bundle $f: X\to Y$.
\item $\chr k>0$ and $K_X$ is nef.
\end{enumerate}
If the characteristic of $k$ is $0$ or larger than $7$ and we are in cases 1-3 then
$X$ is rationally connected.
\end{proposition}

\begin{proof}
Let us recall that if $\chr k=0$ then $X$ is uniruled and $K_X$ is not nef. So we are either in case 4 or $K_X$ is not nef.  By Koll\'ar's and Mori's theorem (see \cite[Main Theorem]{Kol-3}) there exists  a fibration $f: X\to Y$, which is the contraction of a negative extremal ray. By Lemma \ref{contractions} $X$ does not have divisorial contractions. By classification  $f$ is   $f$ is of Fano type (there are no small contractions of smooth $3$-folds). If $\dim Y=0$ we are in the first case. If $\dim Y=1$ then $Y=\PP ^1$ by Theorem \ref{simply-connected}. In this case any fiber of $f$ with reduced scheme structure is a del Pezzo surface. Let us note that all del Pezzo surfaces are rationally connected.
If $\dim Y=2$ then $Y$ is smooth and $f$ is  a flat conic bundle (see \cite[Main Theorem]{Kol-3}).
In this case Theorem \ref{classification-2} allows us to classify possible surfaces  $Y$. In particular, 
If the characteristic of $k$ is $0$ or larger than $7$ then $Y$ is rationally connected. Since a general fiber of $f$ 
is rationally connected, $X$ is also rationally connected.
\end{proof}

\medskip

\begin{remark}
In the first case one knows the classification of such Fano $3$-folds in an arbitrary characteristic (see \cite{SB}).
In the second case it is known by the results of Patakfalvi--Waldron and Fanelli--Schr\"oer that the generic fiber of $f$
is geometrically normal. Unfortunately, these results do not help much in a full classification of smooth projective D-affine $3$-folds. This problem seems to require some new techniques or a non-trivial case by case analysis.
\end{remark}

}

\section*{Acknowledgements}

The author was partially supported by Polish National Centre (NCN) contract numbers 
2015/17/B/ST1/02634 and 2018/29/B/ST1/01232. The author would like to thank P. Achinger, H. Esnault, N. Lauritzen, Y. Lee and D. Rumynin for some useful remarks.

\end{document}